\documentclass[a4paper,12pt]{amsart}
\usepackage[margin=3cm]{geometry}

\newcommand{\dbar}{\bar\partial}
\newcommand{\inp}[2]{\langle #1 , #2 \rangle}

\newcommand{\vnorm}[1]{\left\|  #1 \right\|}

\newcommand{\Cn}{\mathbb{C}^n}

\newcommand{\C}{\mathbb{C}}

\newcommand{\RN}{\mathbb{R}^N}
\newcommand{\R}{\mathbb{R}}

\newcommand{\N}{\mathbb{N}}
\newcommand{\Nn}{\mathbb{N}^n}

\newcommand{\dopt}[2]{\frac{\partial #1}{\partial #2}}

\newlength{\extendaxesby}

\DeclareMathOperator{\spanc}{span}

\DeclareMathOperator{\trace}{trace}

\newtheorem{thm}{Theorem}
\newtheorem{lem}[thm]{Lemma}
\newtheorem{prop}[thm]{Proposition}
\newtheorem{cor}[thm]{Corollary}

\theoremstyle{definition}

\title{Spectral properties of the 
canonical solution operator to $\bar\partial$}
\author{Friedrich Haslinger}
\address{Universit\"at Wien, Fakult\"at f\"ur Mathematik, Nordbergstrasse 15, A-1090 Wien, \"Osterreich}
\email{friedrich.haslinger@univie.ac.at\\lamelb@member.ams.org}
\thanks{The first author was supported by the FWF, Projekt P19147}
\author{Bernhard Lamel}
\thanks{The second author was supported by the FWF, Projekt P17111}
\begin{document}
\maketitle
\begin{abstract}
  We study boundedness, compactness, and Schatten-class membership of 
  the canonical solution operator to $\dbar$, restricted to $(0,1)$-forms with 
  holomorphic coefficients, on $L^2(d\mu)$ where $\mu$ is 
  a measure with the property that the monomials form an orthogonal family in 
  $L^2 (d\mu)$. The characterizations are formulated in terms of moment properties
  of $\mu$. Our results generalize the results of the first author to 
  several variables, contain some known results for several variables, and
  also cover new ground.
\end{abstract}
\section{Introduction and statement of results}
In this paper, we study spectral properties of the canonical solution 
operator to $\dbar$ acting on spaces of $(0,1)$-forms
with holomorphic coefficients in $L^2 (d\mu)$ for measures $\mu$ with
the property that the monomials $z^\alpha$, 
$\alpha\in\Nn$, are orthogonal in $L^2 (d\mu)$.
This situation covers a number of basic examples:

\begin{itemize}
  \item Lebesgue 
    measure on bounded domains in $\Cn$ which 
    are invariant under the torus action
    \[ (\theta_1, \dots, \theta_n) (z_1, \dots, z_n)
    \mapsto (e^{i\theta_1} z_1,\dots e^{i\theta_n} z_n)\]
    (i.e. Reinhardt domains).
  \item Weighted $L^2$ spaces with 
    radially symmetric weights (e.g., generalized Fock spaces). 
  \item Weighted 
    $L^2$ spaces with 
    decoupled radial weights, that is, 
    \[ d\mu = e^{\sum_j \varphi_j (|z_j|^2)} dV, \]
    where $\varphi_j \colon \R \to \R$ is a weight function.
\end{itemize}

Sufficient conditions for the weight in order for the Fock space to be infinite dimensional 
are known from the work of Shigekawa \cite{Shi91}. 
Some of these examples have been studied previously; our approach has
the advantage of unifying these previous result as well as of being
applicable in new situations as well. Our main focus in this paper is 
the case $n > 1$; indeed, we generalize results of the first author
(see \cite{Has02}, \cite{Has01}, \cite{Has01_2}) to this setting.

The behaviour of the canonical solution operator $S$ is interesting
from many points of view. First, there is a close connection between 
properties of $S$ and properties of the $\dbar$-Neumann operator $N$;
indeed, $S = \dbar^* N$. In particular, noncompactness of $S$ prohibits 
compactness of $N$. As is well known, $S$ behaves quite nicely on 
spaces of $(0,1)$-forms with {\em holomorphic} coefficients, and we shall exploit
this connection. On the other hand, for {\em convex} domains, a result of Fu and
Straube \cite{FuStr98} shows that compactness of $S$ on forms with
holomorphic coefficients is also {\em sufficient} for compactness on 
all of $L^2$. 

There is also an intriguing connection between 
the canonical solution operator $S$ and the theory of magnetic 
Schr\"odinger operators (see \cite{FuStr02} and \cite{Has06}); this connection has been exploited
in the recent paper of the first author and Helffer \cite{HaHe06} in order
to study compactness of $S$ on general (not rotation-invariant) weighted
$L^2$-spaces on $\Cn$. 

Let us introduce the notation used in this paper. We denote by
\[ A^{2} (d \mu) = \overline{\left\{ z^\alpha \colon \alpha\in \Nn \right\}},\]
the closure of the monomials in $L^{2} (d\mu)$, 
and write 
\[ m_\alpha = c_{\alpha}^{-1} = \int |z^\alpha|^2 d\mu.\]
We will give necessary and sufficient
conditions in terms of these multimoments of the measure
$\mu$ for the canonical solution operator  to $\bar \partial$, when 
restricted to $(0,1)$-forms with coefficients in $A^{2}(d\mu)$ to be bounded, 
compact, and to belong to the Schatten class $\mathcal{S}_p$. 
This is accomplished by presenting 
a complete diagonalization of the solution operator 
by orthonormal bases with corresponding estimates.
In the case of radially symmetric measures our results specializes to the 
results of \cite{LoY04} applied to this specific case; we are also able to characterize
membership in $\mathcal{S}_p$ for all positive $p$ in some cases 
(a question left open in \cite{LoY04}).

As usual, for a given function space $\mathcal{F}$, $\mathcal{F}_{(0,1)}$ denotes 
the space of $(0,1)$-forms with coefficients in $\mathcal{F}$, that is, expressions
of the form
\[ \sum_{j=0}^n f_j d\bar z_j, \quad f_j \in \mathcal{F}.\]
The $\bar\partial$ operator is the densely defined operator 
\[ \bar\partial f = \sum_{j=1}^n \dopt{f}{\bar z_j} d\bar z_j.\]

The {\em canonical solution operator} $S$ assigns
to each $\omega\in L^2_{(0,1)} (d\mu)$ the solution to the $\bar\partial$ 
equation which is orthogonal to $A^2 (d\mu)$; this solution need not
exist, but if the $\bar\partial$ equation for $\omega$ can be
solved, then $S\omega$ is defined, and is given by 
the unique $f\in L^{2} (d\mu )$
which satisfies
\[ \bar\partial f = \omega \text{ in the sense of distributions and }
f \perp A^{2} (d\mu).\]

Our main interest in this paper is the spectral behaviour of the map $S$ 
restricted to $A^{2}_{(0,1)} (d\mu)$. We first give a criterion for $S$ to be 
a bounded operator. We will frequently encounter multiindeces $\gamma$
which might have
one (but not more than one)
entry equal to $-1$: in that case, we define $c_\gamma = 0$. We will denote the 
set of these multiindeces by $\Gamma$. We let $e_j = (0,\cdots,1,\cdots,0)$ be the 
multiindex with a $1$ in the $j$th spot and $0$ elsewhere.

\begin{thm}\label{t:bounded}
  $S\colon A^{2}_{(0,1)} (d\mu) \to L^{2}(d\mu)$ is bounded if and only if there
  exists a constant $C$ such that 
  \[ \frac{c_{\gamma + e_p}}{c_{\gamma+2e_p} } - \frac{c_{\gamma}}{c_{\gamma+e_p}}
  < C \]
  for all multiindeces $\gamma\in\Gamma$.
\end{thm}

We have a similar criterion for compactness:

\begin{thm}\label{t:compact}
  $S\colon A^{2}_{(0,1)} (d\mu) \to L^{2}(d\mu)$ is compact if and only if 
  \begin{equation} 
    \label{e:limcon}
    \lim_\gamma 
    \left(\frac{c_{\gamma + e_p}}{c_{\gamma+2e_p} } - \frac{c_{\gamma}}{c_{\gamma+e_p}}
  \right)  = 0
  \end{equation}
  for all $p=1,\cdots,n$.
\end{thm}

In particular, the only if implication of 
Theorem~\ref{t:compact} implies several known noncompactness statements for $S$, e.g. 
of Knirsch and Schneider \cite{KnSchn06}, Schneider \cite{Schn06}, as well as 
the noncompactness of $S$ on the polydisc. The main interest in these
noncompactness statements is that if $S$ fails to be compact, so does the $\dbar$-Neumann
operator $N$. 

The multimoments also lend themselves to characterizing the finer spectral 
property of being in the Schatten class $\mathcal{S}_p$. Let us recall that
an operator $T:H_1 \to H_2$ belongs to the Schatten class $\mathcal{S}_p$ if 
the self-adjoint operator $T^{*}T$ has a sequence of eigenvalues belonging to
$\ell^p$. 

\begin{thm}  \label{t:schatten} Let $p>0$. Then $S\colon A^{2}_{(0,1)} (d\mu) \to A^2 (d\mu)$ is in the Schatten-$p$-class
  $\mathcal{S}_p$ if and only if 
  \begin{equation}
    \sum_{\gamma\in\Gamma} \left( \sum_{j}\frac{c_{\gamma+e_j}}{c_{\gamma+2e_j}} - 
    \frac{c_{\gamma}}{c_{\gamma+e_j}}\right)^{\frac{p}{2}} < \infty
    \label{e:schatten}
  \end{equation}
\end{thm}

The condition above is substantially easier to check if $p=2$ (we will show 
that the sum 
is actually a telescoping sum then), i.e. for the case of the 
Hilbert-Schmidt class; we state this as a Theorem: 

\begin{thm}
  \label{t:hilbertschmidt} The canonical solution operator $S$ is in the 
  Hilbert-Schmidt class if and only if
  \begin{equation}
    \lim_{k\to\infty} 
    \sum_{\substack{\gamma\in \Nn, |\gamma|=k \\ 1\leq p \leq n}}  
    \frac{c_{\gamma}}{c_{\gamma+e_p} } < \infty.
    \label{e:hslimcon}
  \end{equation}
\end{thm}

\subsection{Application in the case of decoupled weights}
Let us apply Theorem~\ref{t:bounded} to the case of decoupled weights, or more generally, of product 
measures $d\mu = d\mu_1 \times\cdots\times d\mu_n$, where each $d\mu_j$ is a (circle-invariant) measure on $\C$. Note that 
for such measures, there is {\em definitely} no compactness by Theorem~\ref{t:compact}. If we denote 
by 
\[ c^j_k = \left( \int_\C |z|^{2j} d\mu_k \right)^{-1}, \]
we have that 
\[ c_{(\gamma_1,\cdots,\gamma_n)} = \prod_{k=1}^n c_k^{\gamma_k}. \]
We thus obtain the following corollary. 
\begin{cor} For a product measure $d\mu = d\mu_1\times\cdots \times d\mu_n$ as above, the canonical solution 
  operator
   $S\colon A^{2}_{(0,1)} (d\mu) \to L^{2}(d\mu)$ is bounded if and only if there
  exists a constant $C$ such that 
  \[ \frac{c^{j+1}_k}{c^{j+2}_{k} } - \frac{c^{j}_{k}}{c^{j+1}_{k}}
  < C \]
  for all $j\in \N$ and for all $k=1,\cdots,n$. Equivalently, $S$ is bounded if and only if the 
  canonical solution operator $S_j \colon A^{2} (d\mu_j) \to L^{2}(d\mu_j) $ is bounded for every 
  $j = 1,\cdots,n$.
\end{cor}

\subsection{Application in the case of rotation-invariant measures}

In the case of a rotation-invariant measure $\mu$, we write
\[ m_d = \int_{\Cn} |z|^{2d} d\mu; \] a computation 
(see \cite[Lemma 2.1]{LoY04}) implies that 
\begin{equation}\label{e:cinm}
  c_{\gamma} = \frac{(n + |\gamma| -1)!}{(n-1)!\gamma!} \frac{1}{m_{|\gamma|}}.
\end{equation}
In order to express the conditions of our Theorems, we compute (setting $d= |\gamma| +1$)
\begin{equation}
  \sum_{p} \left( \frac{c_{\gamma+e_p}}{c_{\gamma+2e_p}}  - 
  \frac{c_{\gamma}}{c_{\gamma+e_p}} \right) = 
  \begin{cases}
    \frac{d+2n -1}{d+n} \frac{m_{d+1}}{m_{d}} - \frac{m_{d}}{m_{d-1}} & \gamma_p \neq -1 \text{ for all } p \\
    \frac{1}{d+n} \frac{m_{d+1}}{m_{d}} & \text{else.}
  \end{cases}
  \label{e:expressionmult}
\end{equation}
Note that the Cauchy-Schwarz inequality implies that for large enough $d$, the 
first case in \eqref{e:expressionmult} always dominates the second case; using
this observation and  some trivial inequalities, we get 
the following Corollaries, which should be compared to the results of the first 
author in the one-dimensional case \cite{Has02} and the results of Lovera-Youssfi \cite{LoY04}.

\begin{cor}
  \label{c:boundedrot} Let $\mu$ be a rotation invariant measure on $\Cn$. Then the 
  canonical solution operator to $\dbar$ is bounded on $A^2_{(0,1)}(d\mu)$ if and only if 
  \begin{equation}
    \sup_{d\in\N} 
    \left(
    \frac{(2n+d-1) m_{d+1}}{(n+d) m_d} - \frac{m_d}{m_{d-1}} \right) < \infty
    \label{e:boundedrot}
  \end{equation}
\end{cor}

\begin{cor}
  \label{c:compactrot} 
  Let $\mu$ be a rotation invariant measure on $\Cn$. Then the 
  canonical solution operator to $\dbar$ is compact on $A^2_{(0,1)}(d\mu)$ if and only if 
  \begin{equation}
    \lim_{d\to\infty} 
    \left(
    \frac{(2n+d-1){m_{d+1}}}{(n+d) m_d} - \frac{m_d}{m_{d-1}} \right) = 0.
    \label{e:compactrot}
  \end{equation}
\end{cor}

\begin{cor} 
  \label{c:hilbertrot}
  Let $\mu$ be a rotation invariant measure on $\Cn$. Then the 
  canonical solution operator to $\dbar$ is a Hilbert-Schmidt operator
  on $A^2_{(0,1)}(d\mu)$ if and only if 
  \begin{equation}
    \lim_{d\to\infty} 
    \binom{n+d-1}{n-1} \frac{m_{d+1}}{m_d} < \infty.
    \label{e:hilbertrot}
  \end{equation}
\end{cor}

\begin{cor} 
  \label{c:schattenrot}
  Let $\mu$ be a rotation invariant measure on $\Cn$, $p>0$. Then the 
  canonical solution operator to $\dbar$ is in the Schatten class 
  $\mathcal{S}_p$, as an operator  
  from $A^2_{(0,1)}(d\mu)$ to $L^2 (d\mu)$ if and only if 
  \begin{equation}
    \sum_{d=1}^{\infty} 
    \binom{n+d-2}{n-1} \left(
    \frac{(2n+d-1)m_{d+1}}{(n+d) m_d}  - \frac{m_d}{m_{d-1}} \right)^{\frac{p}{2}} 
    < \infty.
    \label{e:schattenrot}
  \end{equation}
\end{cor}

In particular, Corollary~\ref{c:schattenrot} improves Theorem C of \cite{LoY04} in the sense
that it also covers the case $0< p <2$. We would like to note that our techniques 
can be adapted 
to the setting of \cite{LoY04} by considering the canonical solution operator on 
a Hilbert space $\mathcal{H}$ of holomorphic functions endowed with a norm which is comparable
to the $L^2$-norm on each subspace generated by monomials of a fixed degree $d$, if 
in addition to the requirements in \cite{LoY04} we also assume that
the monomials belong to $\mathcal{H}$; this
introduces the additional weights found by \cite{LoY04} in the formulas, 
as the reader can check. 
In our setting, the formulas are
somewhat ``cleaner'' by working with $A^{2} (d\mu)$ (in particular, Corollary~\ref{c:hilbertrot} only 
holds in this setting). 
\section{Monomial bases and diagonalization}

In what follows, we will denote by 
\[ u_{\alpha} = \sqrt{c_{\alpha}}{z^{\alpha}} \]
the orthonormal basis  of monomials for the space $A^2 (d\mu)$, and by
$U_{\alpha,j} =
u_\alpha d\bar z_j$ the corresponding basis of $A^{2}_{(0,1)} (d\mu)$. 
We first note that it is always possible to solve the $\bar\partial$-equation
for the elements of this basis; 
indeed, $\bar \partial \bar z_j u_\alpha = U_{\alpha,j}$. 
The canonical solution operator is also easily determined for forms
with monomial coefficients:
\begin{lem}
  \label{l:canonicalmonoms} The canonical solution $Sz^{\alpha}d \bar z_j$
  for monomial forms is given by
  \begin{equation}
    S z^{\alpha} d \bar z_j = \bar z_j z^\alpha - 
    \frac{c_{\alpha-e_j}}{c_{\alpha}} z^{\alpha-e_j}, \quad \alpha\in\Nn. 
    \label{e:canonicalmonoms}
  \end{equation}
\end{lem}
\begin{proof}
  We have $\inp{\bar z_j z^\alpha}{z^\beta} = \inp{z^\alpha}{z^{\beta+e_j}}$;
  so this expression is nonzero only if $\beta = \alpha - e_j$ (in particular,
  if this implies \eqref{e:canonicalmonoms} for multiindeces
  $\alpha$ with $\alpha_j =0$; recall our convention that
  $c_{\gamma} = 0 $ if one of the entries of $\gamma$ is negative). Thus
  $S z^{\alpha} d \bar z_j = \bar z_j z^{\alpha} + c z^{\alpha - e_j}$, and
  $c$ is computed by 
  \[ 0 = \inp{\bar z_j z^{\alpha} + c z^{\alpha-e_j}}{z^{\alpha- e_j}}
  = c_{\alpha}^{-1} + c c_{\alpha-e_j}^{-1},\] which 
  gives $c = - c_{\alpha-e_j}/c_{\alpha}$.
\end{proof}
We are going to introduce an orthogonal decomposition
\[ A^{2}_{(0,1)} (d\mu) = \bigoplus_{\gamma\in\Gamma} E_{\gamma} \]
of $A^{2}_{(0,1)}(d\mu)$ into at most $n$-dimensional subspaces
$E_{\gamma}$ indexed by multiindeces $\gamma\in\Gamma$ (we will describe
the index set below), and a corresponding
sequence of mutually orthogonal finite-dimensional subspaces
$F_{\gamma}\subset L^{2}(d\mu)$
which diagonalizes $S$ (by this we mean that $SE_\gamma = F_\gamma$). 
To motivate the definition
of $E_{\gamma}$, note that 
\begin{equation}\label{e:innerprod} \inp{Sz^\alpha d \bar z_k}{S z^{\beta} d\bar z_\ell} = 
\begin{cases}
  0 &  \beta \neq \alpha + e_\ell - e_k, \\
  \frac{1}{c_{\alpha}} \left(  
  \frac{c_{\alpha}}{c_{\alpha+e_\ell}} - \frac{c_{\alpha-e_k}}{c_{\alpha+e_\ell - e_k}}\right)
   &   \beta= \alpha + e_\ell - e_k,
\end{cases}
\end{equation}
so that $\inp{Sz^\alpha d \bar z_k}{S z^{\beta} d\bar z_\ell} \neq 0$ if and 
only if there exists a multiindex $\gamma $ such that $\alpha = \gamma + e_k$ and
$\beta = \gamma + e_\ell$. We thus define
\[ E_{\gamma} = \spanc \left\{ 
U_{\gamma+e_j,j} \colon 1\leq j\leq n\right\}=
\spanc \left\{ 
z^{\gamma+e_j} d\bar z_j \colon 1\leq j\leq n\right\},\]
and likewise $F_\gamma = S E_\gamma$. Recall that
$\Gamma$ is defined 
to be the set of all multiindeces
whose entries are greater or equal to $-1$ and at
most one negative entry. Note that
$E_{\gamma}$ is $1$-dimensional
if exactly one entry in $\gamma$ equals $-1$, 
and $n$-dimensional otherwise.
We have already observed that $F_\gamma$ are mutually orthogonal subspaces
of $L^{2} (d\mu)$. 

Whenever we use multiindeces $\gamma$ and
integers $p\in\left\{ 1,\cdots,n \right\}$ as indeces, 
we use the convention that the $p$ run over all $p$ such 
that $\gamma+e_p \geq 0$; that is, for a fixed multiindex $\gamma\in\Gamma$, 
either the indeces are either
 all $p\in \left\{ 1,\cdots,n \right\}$ or there
is exactly one $p$ such that $\gamma_p = -1$, in which case the index
is exactly this one $p$.

We next observe that we can find an orthonormal basis of $E_{\gamma}$ and an 
orthonormal basis of $F_\gamma$ such that in these bases
$S_\gamma = S|_{E_\gamma}\colon E_\gamma \to F_\gamma$ acts diagonally. 
First note that it is enough to do this if $\dim E_\gamma = n$ (since an operator
between one-dimensional spaces is automatically diagonal).
Fixing $\gamma$, the functions $U_j := U_{\gamma+e_j,j}$ are an orthonormal basis of $E_\gamma$. The 
operator $S_\gamma$ is clearly nonsingular on this space, so 
the functions $S U_{j}=\Psi_j$ constitute
a basis of $F_\gamma$. For a basis $B$ of vectors $v^j = 
\left( v^j_1,\dots,v^j_n \right)$, $j = 1,\dots,n$ of $\Cn$ we consider
the new basis 
\[ V_k = \sum_{j=1}^{n} v_k^j U_{j};\]
since the basis given by the $U_{j}$ is orthonormal, 
the basis given by the $V_k$ is also orthonormal provided that
the vectors $v_k = (v_k^1,\cdots,v_k^n)$ constitute an orthonormal
basis for $\Cn$ with the standard hermitian product.
Let us write 
\[ \Phi_k = S V_k = \sum_j v^j_k S U_{j}. \]
The inner product $\inp{\Phi_p}{\Phi_q}$ is then given by 
$\sum_{j,k} v_p^j \bar v_q^k \inp{S U_{j}}{S U_{k}}$. 
We therefore have 
\begin{multline}
\begin{pmatrix}
  \inp{\Phi_1}{\Phi_1} & \cdots & \inp{\Phi_1}{\Phi_n} \\
  \vdots & & \vdots \\
  \inp{\Phi_n}{\Phi_1} & \cdots & \inp{\Phi_n}{\Phi_n}
\end{pmatrix}=\\
\begin{pmatrix}
  v_1^1 & \cdots & v_1^n \\
  \vdots & & \vdots \\
  v_n^1 & \cdots & v_n^n
\end{pmatrix}
\begin{pmatrix}
  \inp{\Psi_1}{\Psi_1} & \cdots & \inp{\Psi_1}{\Psi_n} \\
  \vdots & & \vdots \\
  \inp{\Psi_n}{\Psi_1} & \cdots & \inp{\Psi_n}{\Psi_n}
\end{pmatrix}
\begin{pmatrix}
  \bar v_1^1 & \cdots & \bar v_n^1 \\
  \vdots & & \vdots \\
  \bar v_1^n & & \bar v_n^n
\end{pmatrix}.
\end{multline}
Since the matrix $(\inp{\Psi_j}{\Psi_k})_{j,k}$ is hermitian, 
we can unitarily diagonalize it; that is, we can choose an orthnormal basis $B$ of $\Cn$ such 
that with this
choice of $B$ the vectors $\varphi_{\gamma,k} =
V_k = \sum_j v^j_k U_{\gamma+e_j,j}$ of $E_\gamma$ are orthonormal, and their 
images $\Phi_{k}= S V_k$ are orthogonal in $F_{\gamma}$.
Therefore, 
$\Phi_k / \| \Phi_k \|$ is an orthonormal basis of $F_\gamma$ such that 
$S_\gamma: E_\gamma\to F_\gamma$ is diagonal when
expressed in terms of the bases 
$\left\{ V_1,\cdots,V_n \right\} 
\subset E_{\gamma}$ and $\left\{ \Phi_1,\cdots,\Phi_n \right\}\subset F_{\gamma}$, 
with entries $\| \Phi_k \|$. 

Furthermore, the $\| \Phi_k \|$ are exactly the square roots of 
the eigenvalues of the matrix
$(\inp{\Psi_p}{\Psi_q})$ which by \eqref{e:innerprod} is given by 
\begin{equation}
  \begin{aligned}
  \inp{\Psi_p}{\Psi_q} &= 
  \inp{SU_{\gamma+e_p,p}}{SU_{\gamma+e_q,q}} \\
  &= \sqrt{c_{\gamma+e_p}}\sqrt{c_{\gamma+e_q}} \inp{S\, z^{\gamma+e_p}\,d\bar z_p}{S\, z^{\gamma+e_q} \, d\bar z_q} \\
  &=  \sqrt{c_{\gamma+e_p}c_{\gamma+e_q}} \frac{1}{c_{\gamma+e_p}}
  \left( \frac{c_{\gamma+e_p}}{c_{\gamma+e_p +e_q}} - \frac{c_{\gamma}}{c_{\gamma+e_q}}  \right) \\
  &=
  \frac{c_{\gamma+e_p}c_{\gamma+e_q} - c_{\gamma}
  c_{\gamma+e_p+e_q}}{c_{\gamma+e_p + e_q}\sqrt{c_{\gamma+e_p}c_{\gamma+e_q}}}
\end{aligned}
\end{equation}

Summarizing, we have the following Proposition.
\begin{prop}\label{p:diagonal}
  With $\mu$ as above, the canonical solution operator $S\colon
  A^{2}_{(0,1)} (d\mu) \to L^2_{(0,1)} (d\mu)$ admits a diagonalization by orthonormal
  bases. In fact, we have a decomposition $A^{2}_{(0,1)} = \bigoplus_{\gamma} E_\gamma$
  into mutually orthogonal finite dimensional
  subspaces $E_\gamma$, indexed by the multiindeces $\gamma$ 
  with at most one negative entry (equal to $-1$), which are of dimension $1$ or $n$,
  and orthonormal bases $\varphi_{\gamma,j}$
  of $E_{\gamma}$,
  such that $S \varphi_{\gamma,j} $ is a set of mutually orthogonal vectors
  in $L^{2} (d\mu)$. For fixed $\gamma$, the 
  norms $\vnorm{ S_{} \varphi_{\gamma,j}}$ are the 
  square roots of the eigenvalues of the matrix $C_{\gamma} = (C_{\gamma,p,q})_{p,q}$
  given by 
  \begin{equation}
    C_{\gamma,p,q} =  \frac{c_{\gamma+e_p}c_{\gamma+e_q} - c_{\gamma}
  c_{\gamma+e_p+e_q}}{c_{\gamma+e_p + e_q}\sqrt{c_{\gamma+e_p}c_{\gamma+e_q}}}.
    \label{e:inpmat2}
  \end{equation}
  In particular, we have that 
  \begin{equation}
    \sum_{j=1}^{n} \vnorm{S \varphi_{\gamma,j}}^2 =
    \trace (C_{\gamma,p,q})_{p,q}
=    \sum_{p=1}^{n} \left( \frac{c_{\gamma+e_p}}{c_{\gamma+2e_p}} - \frac{c_{\gamma}}{c_{\gamma+e_p} } \right)
    \label{e:traces}
  \end{equation}
\end{prop}
\section{Boundedness: Proof of Theorem~\ref{t:bounded}}
In order to prove Theorem~\ref{t:bounded}, we are using 
Proposition~\ref{p:diagonal}. We have seen that we have an orthonormal basis
$\varphi_{\gamma,j}$, $\gamma\in\Gamma$, $j\in\left\{ 1,\cdots,
n\right\}$, such 
that the images 
$S \varphi_{\gamma,j}$ are mutually orthogonal. Thus, $S$ is bounded
if and only if there exists a constant $C$ such that
\[ \vnorm{S\varphi_{\gamma,j}}^2 \leq C \] 
for all $\gamma\in\Gamma$ and $j\in\left\{ 1,\cdots,\dim E_{\gamma}\right\}$.
If $\dim E_\gamma = 1$, then $\gamma$ has exactly one entry (say the $j$th one)
equal to $-1$; in that case, let us write $\varphi_\gamma = 
U_{\gamma+e_j} d\bar z_j$. We have 
$S\varphi_\gamma = \sqrt{c_{\gamma+e_j}}\bar z_j z^{\gamma+e_j}$, and so  
\[ \vnorm{S\varphi_\gamma}^2 = \frac{c_{\gamma+e_j}}{c_{\gamma+2e_j}}.\]
On the other hand, if $\dim E_\gamma = n$, we argue as follows:
Writing $\vnorm{S\varphi_{\gamma,j}}^2 = \lambda_{\gamma,j}^2$ with $\lambda_{\gamma,j} > 0$, 
from \eqref{e:traces} we find that
\[ \sum_{j=1}^n \lambda_{\gamma,j}^2 = \sum_{j=1}^n
\left(\frac{c_{\gamma+e_j}}{c_{\gamma+2e_j}} -
\frac{c_{\gamma}}{c_{\gamma+e_j}}\right).\]
The last 2 equations complete the proof of Theorem~\ref{t:bounded}.

\section{Compactness}

In order to prove Theorem~\ref{t:compact}, we use the following elementary Lemma
(which is for example contained in  \cite{DAcarus}):
\begin{lem}
  \label{l:compcrit} Let $H_1$ and $H_2$ be Hilbert spaces, and 
  assume that $S\colon H_1 \to H_2$ is a bounded linear operator. 
  Then $S$ is compact if and only if for every $\varepsilon>0$ there
  exists a compact operator $T_\varepsilon\colon H_1 \to H_2$ such that the
  following
  inequality holds:
  \begin{equation}
    \label{e:compinequ}
    \vnorm{Sv}_{H_2}^2 \leq 
    \vnorm{T_\varepsilon v}_{H_2}^2 + \varepsilon \vnorm{v}_{H_1}^2.
  \end{equation}
\end{lem}

\begin{proof}[Proof of Theorem~\ref{t:compact}]
  We first show that \eqref{e:limcon} implies compactness. 
  We will use the notation which was already used in the proof
  of Theorem~\ref{t:bounded}; that is, we write 
  $\vnorm{S \varphi_{\gamma,j}}^2 = \lambda_{\gamma,j}^2$.
  Let $\varepsilon > 0$. 
  There exists a finite set $A_\varepsilon$ of 
  multiindeces $\gamma\in\Gamma$ such that for all $\gamma
  \notin A_\varepsilon$, 
  \[  
  \sum_{j=1}^n \lambda_{\gamma,j}^2 = \sum_{j=1}^n
  \left(\frac{c_{\gamma+e_j}}{c_{\gamma+2e_j}} -
  \frac{c_{\gamma}}{c_{\gamma+e_j}}\right) < \varepsilon.
  \]
  Hence, if we consider the finite dimensional (and thus, compact) operator
  $T_{\varepsilon}$ defined by 
  \[ T_{\varepsilon} \sum a_{\gamma,j} \varphi_{\gamma,j} =
  \sum_{\gamma\in A_{\varepsilon}} a_{\gamma,j} S \varphi_{\gamma,j}, \]
  for any $v = \sum a_{\gamma,j} \varphi_{\gamma,j}\in A^2_{(0,1)} (d\mu)$
  we obtain
  \[\begin{aligned} 
    \vnorm{S v}^2 
    &= 
    \vnorm{T_{\varepsilon}  v}^2
    + \vnorm{ 
    S\sum_{\gamma\notin A_{\varepsilon}} a_{\gamma,j} \varphi_{\gamma,j}}^2 \\
    &= \vnorm{T_{\varepsilon} v}^2
    + \sum_{\gamma\notin A_{\varepsilon}} 
    |a_{\gamma,j}|^2 \vnorm{S\varphi_{\gamma,j}}^2 
    \\
    &= \vnorm{T_{\varepsilon} v}^2 + 
    \sum_{\gamma\notin A_{\varepsilon}}
    |a_{\gamma,j}|^2 \lambda_{\gamma,j}^2 \\
    &\leq \vnorm{T_{\varepsilon} v}^2 + \varepsilon
    \sum_{\gamma\notin A_{\varepsilon}} |a_{\gamma,j}|^2 \\
    & \leq \vnorm{T_{\varepsilon} v}^2 + \varepsilon \vnorm{v}^2.
  \end{aligned}\]
  Hence, \eqref{e:compinequ} holds and we have proved the first implication
  in Theorem~\ref{t:compact}.

  We now turn to the other direction. Assume that \eqref{e:limcon} is not 
  satisfied. Then there exists a $K>0$ and an infinite family $A$ of multiindeces
  $\gamma$ such that for all $\gamma\in A$,
  \[  
  \sum_{j=1}^n \lambda_{\gamma,j}^2 = \sum_{j=1}^n
  \left(\frac{c_{\gamma+e_j}}{c_{\gamma+2e_j}} -
  \frac{c_{\gamma}}{c_{\gamma+e_j}}\right)> nK.
  \]
  In particular, for each $\gamma \in A$, there exists a $j_{\gamma}$ such 
  that $\lambda_{\gamma,j_{\gamma}}^2 > K$. Thus, we have an 
  infinite orthonormal 
  family $\{ \varphi_{\gamma,j_{\gamma}} \colon \gamma\in A \}$ 
  of vectors such that their images $S \varphi_{\gamma,j_{\gamma}}$ are
  orthogonal and
  have norm bounded from below by $\sqrt{K}$, which contradicts compactness.
\end{proof}

\section{Membership in the Schatten classes $\mathcal{S}_p$ and in the Hilbert-Schmidt class}

We keep the notation introduced in the previous sections. 
We will also need to introduce the usual grading on the index set $\Gamma$, that is, we write
\begin{equation}
  \Gamma_k = \left\{ \gamma\in\Gamma \colon |\gamma| = k \right\}, \quad
  k \geq -1.
  \label{e:gammakdef}
\end{equation}

In order to study the membership in the Schatten class, we need the following elementary Lemma:
\begin{lem}\label{l:comparable}
  Assume that $p(x)$ and $q(x)$ are continuous, real-valued functions on $\R^N$ which are homogeneous of degree
  $1$ (i.e. $p(tx) = t p(x)$  and $q(tx) = tq(x)$ for  $t\in\R$), and $q(x) = 0$ as well as $p(x)=0$ 
  implies $x=0$. Then there exists a constant $C$ such that 
  \begin{equation}
    \frac{1}{C} | q(x) | \leq |p(x)| \leq C | q(x) | .
    \label{e:comparable1}
  \end{equation}
\end{lem}
\begin{proof}
  Note that the set $B_q=\left\{ x \colon q(x) = 1 \right\}$ is compact: it's closed since $q$ is continuous, and since $|q|$ is bounded from 
  below on $S^{N}$ by some $m>0$, it is necessarily contained in the closed ball of radius $1/m$. Now, the function $|p|$ is 
  bounded on the compact set $B_q$; say, by $1/C$ from below and $C$ from above. Thus for all $x\in\RN$, 
  \[\frac{1}{C} \leq \left| p\left( \frac{x}{q(x)} \right)\right| \leq C,\]
  which proves \eqref{e:comparable1}.
\end{proof}
\begin{proof}[Proof of Theorem~\ref{t:schatten}]
  Note that $S$ is in the Schatten class $\mathcal{S}_p$ if and only if 
  \begin{equation}
    \sum_{{\gamma\in\Gamma,\,j}} \lambda_{\gamma,j}^{p} < \infty.
    \label{e:inschattenclass}
  \end{equation}
  We rewrite this sum as 
  \[ \sum_{\gamma\in\Gamma} \left( \sum_{j} \lambda_{\gamma,j}^{p} \right) =: M \in \R\cup\left\{ \infty \right\}. \]
  Lemma~\ref{l:comparable}  implies that there exists a constant $C$ such that for every $\gamma\in\Gamma$, 
  \[ \frac{1}{C} \left(\sum_{j} \lambda_{\gamma,j}^{2} \right)^{p/2} \leq 
  \sum_{j} \lambda_{\gamma,j}^{p} \leq C \left(\sum_{j} \lambda_{\gamma,j}^{2} \right)^{p/2}. \]
  Hence, $M<\infty$ if and only if 
  \[\sum_{\gamma}  \left(\sum_{j} \lambda_{\gamma,j}^{2} \right)^{p/2} < \infty,\]
  which after applying \eqref{e:traces} becomes the condition \eqref{e:schatten} claimed in Theorem~\ref{t:schatten}.
\end{proof}

\begin{proof}[Proof of Theorem~\ref{t:hilbertschmidt}]
  $S$ is in the 
  Hilbert-Schmidt class if and only if 
  \begin{equation}
    \sum_{\gamma\in\Gamma, j} \lambda_{\gamma,j}^{2} < \infty.
    \label{e:hsorig}
  \end{equation}
  We will prove that 
  \begin{equation}
    \sum_{\ell =-1}^{k} \sum_{\gamma\in\Gamma_\ell, j} \lambda_{\gamma,j}^2 = 
    \sum_{\substack{\alpha\in \Nn, |\alpha|=k+1 \\ 1\leq p \leq n}}  
    \frac{c_{\alpha}}{c_{\alpha+e_p} },
    \label{e:hsequiv}
  \end{equation}
  which immediately implies Theorem~\ref{t:hilbertschmidt}. The proof is by 
  induction over $k$. For $k= - 1$,  the left hand side of \eqref{e:hsequiv}
  is \[ \sum_{j=1}^n \lambda_{-e_j,j}^2  =  \sum_{j=1}^n 
  \vnorm{z_j}^2 c_0 = \sum_{j=1}^n \frac{c_0}{c_{e_p}},\]
  which is equal to the right hand side.
  Now assume that the \eqref{e:hsequiv} holds for $k = K-1$; we will show that
  this implies it holds for $k = K$. We write
  \[ \begin{aligned}
    \sum_{\ell =-1}^{K} \sum_{\gamma\in\Gamma_\ell, j} \lambda_{\gamma,j}^2 
    &=  \sum_{\substack{\alpha\in \Nn, |\alpha|=K-1 \\ 1\leq p \leq n}}  
    \frac{c_{\alpha}}{c_{\alpha+e_p} } +
    \sum_{\gamma\in\Gamma_{K},j} 
    \left(\frac{c_{\gamma+e_j}}{c_{\gamma+2e_j}} -
    \frac{c_{\gamma}}{c_{\gamma+e_j}}\right) \\
    &= \sum_{\substack{\alpha\in \Nn, |\alpha|=K \\ 1\leq p \leq n}}  
    \frac{c_{\alpha}}{c_{\alpha+e_p} } .
  \end{aligned}\]
  This finishes the proof of Theorem~\ref{t:hilbertschmidt}.
\end{proof}
\bibliographystyle{abbrv}
\bibliography{bibliography}
\end{document}